\documentclass[12pt, twoside, leqno]{article}

\usepackage{amsmath,amsthm}
\usepackage{amssymb}

\usepackage{enumerate}

\usepackage{graphicx}

\usepackage{pgf,tikz}
\usetikzlibrary{arrows}

\newtheorem{thm}{Theorem}[section]

\newtheorem{lem}[thm]{Lemma}

\theoremstyle{definition}
\newtheorem{defin}[thm]{Definition}

\numberwithin{equation}{section}

\usepackage{pgf,tikz}
\usetikzlibrary{arrows}

\begin{document}

\baselineskip=12.1pt

%%%%%%%%%%%%%%%%

\title{ The automorphism group of the Andr$\acute{a}$sfai  graph}

\author{S.Morteza Mirafzal\\
Department of Mathematics \\
  Lorestan University, Khorramabad, Iran\\
E-mail: mirafzal.m@lu.ac.ir\\
E-mail: smortezamirafzal@yahoo.com}

\date{}

\maketitle

\renewcommand{\thefootnote}{}


\footnote{\emph{Date}:  }

\renewcommand{\thefootnote}{\arabic{footnote}}
\setcounter{footnote}{0}
\date{}

\begin{abstract}

Let $k \geq 1$ be an integer and $n=3k-1$.   Let
$\mathbb{Z}_n$ denote the additive group of integers modulo $n$ and let $C$ be
the subset of $\mathbb{Z}_n$ consisting of the elements congruent to 1 modulo 3.
 The Cayley graph $Cay(\mathbb{Z}_n; C)$ is known as the Andr$\acute{a}$sfai graph And($k$). In this note, we  determine the automorphism group of this graph.  We will show that $Aut(And(k))$ is isomorphic with the dihedral group $\mathbb{D}_{2n}$.

\end{abstract}

\maketitle

\section{ Introduction}
 In this paper, a graph $\Gamma=(V,E)$ is
considered as an undirected simple graph where $V=V(\Gamma)$ is the vertex-set
and $E=E(\Gamma)$ is the edge-set. For all the terminology and notation
not defined here, we follow [4,15].\

 Let $m>0$ be an integer.  Let
$\mathbb{Z}_m$ denote the additive group of integers modulo $m$. Let $ k >  1$ be an integer and $n=3k-1$. Let $C=\{3t+1\  | \  0\leq t \leq k-1  \}$ be
the subset of $\mathbb{Z}_n$ consisting of the elements congruent to 1 modulo 3. It is easy to see that $C$ is a symmetric set, that is, $C$ is an inverse closed subset of the group $\mathbb{Z}_n$.
 The Cayley graph $Cay(\mathbb{Z}_n; C)$ is known as the Andr$\acute{a}$sfai graph And($k$).
It is easy to check that the
graph And(2)
is isomorphic to the 5-cycle and the graph And(3) is  the M$\ddot{o}$bius ladder of order 8.  The graph
And(4) is depicted in Figure 1.

\definecolor{qqqqff}{rgb}{0.,0.,1.}
\begin{tikzpicture}[line cap=round,line join=round,>=triangle 45,x=0.90cm,y=0.8cm]
\clip(-3.42,-1.3) rectangle (12.2,7.38);
\draw (4.74,3.84)-- (4.82,2.3);
\draw (4.82,2.3)-- (4.28,0.98);
\draw (4.28,0.98)-- (3.18,0.32);
\draw (3.18,0.32)-- (1.76,0.7);
\draw (1.76,0.7)-- (0.66,1.54);
\draw (0.66,1.54)-- (0.06,3.12);
\draw (0.06,3.12)-- (0.26,4.58);
\draw (0.26,4.58)-- (1.06,5.48);
\draw (1.06,5.48)-- (2.52,6.36);
\draw (4.02,5.18)-- (2.52,6.36);
\draw (4.02,5.18)-- (4.74,3.84);
\draw (4.06,5.64) node[anchor=north west] {1};
\draw (5.1,4.12) node[anchor=north west] {2};
\draw (5.2,2.46) node[anchor=north west] {3};
\draw (4.72,0.98) node[anchor=north west] {4};
\draw (3.1,-0.01) node[anchor=north west] {5};
\draw (1.52,0.36) node[anchor=north west] {6};
\draw (0.16,1.62) node[anchor=north west] {7};
\draw (-0.4,3.34) node[anchor=north west] {8};
\draw (-0.32,4.98) node[anchor=north west] {9};
\draw (0.26,5.9) node[anchor=north west] {10};
\draw (4.28,0.98)-- (2.52,6.36);
\draw (2.52,6.36)-- (0.66,1.54);
\draw (4.02,5.18)-- (3.18,0.32);
\draw (4.02,5.18)-- (0.06,3.12);
\draw (4.74,3.84)-- (1.76,0.7);
\draw (4.74,3.84)-- (0.26,4.58);
\draw (4.82,2.3)-- (0.66,1.54);
\draw (4.82,2.3)-- (1.06,5.48);
\draw (4.28,0.98)-- (0.06,3.12);
\draw (3.18,0.32)-- (0.26,4.58);
\draw (1.76,0.7)-- (1.06,5.48);
\draw (2.52,6.92) node[anchor=north west] {0};
\draw (.99,-0.5) node[anchor=north west] {Figure 1. And(4)};
\begin{scriptsize}
\draw [fill=qqqqff] (2.52,6.36) circle (1.5pt);
\draw [fill=qqqqff] (4.74,3.84) circle (1.5pt);
\draw [fill=qqqqff] (4.82,2.3) circle (1.5pt);
\draw [fill=qqqqff] (4.28,0.98) circle (1.5pt);
\draw [fill=qqqqff] (3.18,0.32) circle (1.5pt);
\draw [fill=qqqqff] (1.76,0.7) circle (1.5pt);
\draw [fill=qqqqff] (0.66,1.54) circle (1.5pt);
\draw [fill=qqqqff] (0.06,3.12) circle (1.5pt);
\draw [fill=qqqqff] (0.26,4.58) circle (1.5pt);
\draw [fill=qqqqff] (1.06,5.48) circle (1.5pt);
\draw [fill=qqqqff] (4.02,5.18) circle (1.5pt);
\end{scriptsize}
\end{tikzpicture} \

The Cayley graphs And(k) were first used by Andr$\acute{a}$sfai in [1], and also
appeared in his book [2]. It is not hard to show  that the  graph And(k) has diameter 2 and girth 4.
The Andr$\acute{a}$sfai graph And($k$) has some interesting properties and is a classic  example in the subject of graph homomorphism \cite{4}.  For a given graph, one of the problems concerning   it is  determination   its automorphism group. To the best of our knowledge, the automorphism group of the graph And($k$) is still unknown. The main aim of  the present  paper  is to determine the automorphism group of And($k$). We will show that $Aut(And(k)) \cong \mathbb{D}_{2n}$, $n=3k-1$,  where $\mathbb{D}_{2n}$ denotes the dihedral group of order $2n$.

\section{Preliminaries }

The graphs $\Gamma_1 = (V_1,E_1)$ and $\Gamma_2 =
(V_2,E_2)$ are called $isomorphic$, if there is a bijection $\alpha
: V_1 \longrightarrow V_2 $   such that  $\{a,b\} \in E_1$ if and
only if $\{\alpha(a),\alpha(b)\} \in E_2$ for all $a,b \in V_1$.
In such a case,  the bijection $\alpha$ is called an $isomorphism$.
An $automorphism$ of a graph $\Gamma$ is an isomorphism of $\Gamma
$ with itself. The set of automorphisms of $\Gamma$  with the
operation of composition of functions is a group  called the
$automorphism\  group$ of $\Gamma$ and denoted by $ Aut(\Gamma)$.

 The
group of all permutations of a set $V$ is denoted by $Sym(V)$  or
just $Sym(n)$ when $|V| =n $. A $permutation$ $group$ $G$ on
$V$ is a subgroup of $Sym(V).$  In this case we say that $G$ $acts$
on $V$. If $G$ acts on $V$ we say that $G$ is
$transitive$ on $V$ (or $G$ acts $transitively$ on $V$), when there is just
one orbit. This means that given any two elements $u$ and $v$ of
$V$, there is an element $ \beta $ of  $G$ such that  $\beta (u)= v
$.  If $\Gamma$ is a graph with vertex-set $V$  then we can view
each automorphism of $\Gamma$ as a permutation on $V$  and so $Aut(\Gamma) = G$ is a
permutation group on $V$.

A graph $\Gamma$ is called $vertex$-$transitive$ if  $Aut(\Gamma)$
acts transitively on $V(\Gamma)$. For $v\in V (\Gamma )$ and $G = Aut(\Gamma )$  the $stabilizer\  subgroup$  $G_{v}$ is the subgroup of $G$ containing of all automorphisms  fixing $v$. In the vertex-transitive case all stabilizer subgroups $G_{v}$ are conjugate in $G$, and consequently isomorphic. In this case, the index of $G_{v}$ in $G$ is given by the equation,  $|G : G_{v}| = \frac{|G|}{|G_{v}}| = |V (\Gamma)|$. This fact is known as the Orbit-Stabilizer theorem which is a useful tool in finding the automorphism group of vertex-transitive graphs.\

Let $G$ be any abstract finite group with identity $1$, and
suppose $\Omega$ is a set of   $G$, with the
properties:

(i) $x\in \Omega \Longrightarrow x^{-1} \in \Omega$;  $ \ (ii)
 \ 1\notin \Omega $.

The $Cayley\  graph$  $\Gamma=\Gamma (G; \Omega )$ is the (simple)
graph whose vertex-set and edge-set are defined as follows :

$V(\Gamma) = G $,  $  E(\Gamma)=\{\{g,h\}\mid g^{-1}h\in \Omega \}$.\

 Although in most situations it is difficult to determine the automorphism group
of a graph $G$, and how it acts on its vertex-set or edge-set,   there are various papers in the literature 
and some of the recent works include   [3,5,6,7,8,9,10,11,12,13,14,16,17].

The group $G$ is called a semidirect product of $ N $ by $Q$,
denoted by $ G=N \rtimes Q $,
 if $G$ contains subgroups $ N $ and $ Q $ such that:  (i)$
N \unlhd G $ ($N$ is a normal subgroup of $G$); (ii) $ NQ = G $; and
(iii) $N \cap Q =1 $. \

\section{Main Results}

\begin{defin} Let $ k > 1$ be an integer and $n=3k-1$. Let $C=\{3t+1\  | \  0\leq t \leq k-1  \}$ be
the subset of $\mathbb{Z}_n$ consisting of the elements congruent to 1 modulo 3. It is easy to see that $C$ is a symmetric set, that is, $C$ is an inverse closed subset of the group $\mathbb{Z}_n$.
 The Cayley graph $Cay(\mathbb{Z}_n; C)$ is known as the Andr$\acute{a}$sfai graph And($k$).
\end{defin}

It follows from Definition 3.1, that the graph And(k) is a regular graph of valency $k$ and the vertex-set of And(k)
is the set $V=V_0 \cup V_1 \cup V_2$, where $V_0= \{3t\  | \  0\leq t \leq k-1  \}$, $V_1=\{3t+1\  | \  0\leq t \leq k-1  \}$ and $V_2=\{3t+2\  | \  0\leq t \leq k-2  \}$.
Thus,  we have $|V_0|=|V_1|=k$ and $|V_2|=k-1$. If $v \in
V_0$, then $v=3j$, for some $j$,  $0 \leq j \leq k-1$. Now it is easy to
see that 
$$N(v)=\{3i+1 \ | \ j \leq i \leq k-1 \} \cup \{3l+2 \ | \ 0 \leq l \leq j-1 \}. \ \ \ \ (*)$$

Also, if $w \in V_2$, then $w=3j+2$, $0 \leq j \leq k-2$, and thus we have

$$N(w)=\{3i+1 \ | \ 0 \leq i \leq j \} \cup \{3l \ | \ j+1 \leq l \leq k-1 \}. \ \ \ \  \ \ \ \  (**)$$

Now, from (*) and (**), it follows that the graph induced by the set $V_0 \cup V_2$ in And(k) is a bipartite graph such that the vertex $3j=v \in V_0$ has $j$ neighbors in $V_2$ and the vertex $3j+2=w \in V_2$ has $k-j-1$ neighbors  in $V_0$. Note that all the neighbors of the vertex $v=0$ are in $V_1$.  Let $H=\langle (V_0-\{ 0\}) \cup V_2 \rangle$ be the subgraph induced by the set $(V_0-\{ 0\}) \cup V_2$ in And(k). Thus, $H$ is a connected bipartite graph such that if $v,w$ are distinct vertices in $H$, then we have $N(v) \neq N(w)$ (note that the vertex $v=3(k-1)$ is adjacent to every vertex in $V_2$ and the vertex $w=2$ is adjacent to any vertex in $V_0$).
\\

In the sequel, we need the following fact.

\begin{lem}
Let $\Gamma= (U \cup W,E)$, $U \cap W=\emptyset $ be a connected bipartite graph.  If $f$ is an automorphism of the graph $\Gamma$,  then $ f(U)=U$ and $f(W) =W$,  or $ f(U) = W $ and $f(W) = U$.

\end{lem}

\begin{proof}
Automorphisms of $\Gamma$ preserve
distance between vertices and since two vertices are in the same part if and
only if they are at even distance from each other, the result follows.

\end{proof}

\begin{thm}
Let $k > 1$ be an integer and $n=3k-1$. Then for the automorphism group of the graph And(k)  we have,
$Aut(And(k)) \cong \mathbb{D}_{2n}$,  where $\mathbb{D}_{2n}$ denotes the dihedral group of order $2n$.
\end{thm}

\begin{proof} Let $\Gamma=(V,E)$=And($k$) and $A=Aut(\Gamma)$ be the automorphism group of $\Gamma$. Consider the vertex $v=0$   and let $A_0$ be its stabilizer subgroup, that is, $A_0=\{a\in A \ | \ a(0)=0  \}$. We know that $\Gamma$ is a Cayley graph, hence it is  a vertex-transitive graph.   From the well known Orbit-Stabilizer theorem,  we know    that $|V|=\frac{|A|}{|A_v|}$ and hence $|A|=|V||A_v|$, where $v$ is a vertex in $\Gamma$.  In the first step of our work we determine $|A_0|.$  Let $f\in A_0$.  Let  $V_0$, $V_1$ and $V_2$ be the subsets of $V$ which are defined preceding (*) and $W_0=V_0-\{ 0\}$.   Thus for the restriction of $f$ to $N(0)=V_1$ we have $f(V_1)=V_1$ and hence $f(W_0 \cup V_2)=W_0 \cup V_2$. Let $H$ be the subgraph induced by the set $W_0 \cup V_2$ in the graph $\Gamma=And(k)$ and $g=f|_{W_0 \cup V_2}$. Hence $g$ is an automorphism of $H$. It is clear that $H$ is a connected bipartite graph with parts $W_0$ and $V_2$ such that $|W_0|=|V_2|=k-1$. In each part of the graph $H$ there is exactly one vertex $x_j$ of degree $j$, $1 \leq j \leq k-1$. In other words, the vertex $v_j=3j$ is the unique vertex in $W_0$ of degree $j$, also  the vertex $w_j=3k-1-3j=3(k-j)-1=3(k-j-1)+2$ is the unique vertex in $V_2$ of degree $j$. Note that $w_j=3k-1-3j$ is the inverse of $v_j=3j$ in the cyclic group $\mathbb{Z}_{3k-1}$, hence we   denote it by $-v_j$.
We know that the mapping $g$ is an automorphism of the connected  bipartite graph $H$. Thus from Lemma 3.2, it follows that \newline (i) $g(W_0)=W_0$ or
(ii) $g(W_0)=V_2$. \newline
(i) If $g(W_0)=W_0$, then since the vertex $v=3j$ is the unique vertex of $W_0$ of degree $j$, hence for every $w \in W_0$ we have
$g(w)=w$. Similarly, for  every $v \in V_2$ we have $g(v)=v$. In other words, the restriction of the automorphism $f$ to the set $W_0 \cup V_2$ is the identity mapping. We show that if $x \in V_1$, then $f(x)=x$. Note that we have $f(V_1)=V_1$.  If $v=3j+1$ is a vertex in $V_1$, then the set of neighbors of $v$ in $W_0$ is $N_1= \{ 3t \ | \  \  1\leq t \leq j  \}$. Hence $v$ has exactly $j$ neighbors in $W_0$. Since the number of neighbors of $v$ and $f(v)$ in $W_0$ are equal, hence we must have $v=f(v)$.  From our discussion it follows that if $g(W_0)=W_0$, then the automorphism $f$ is the identity automorphism of the graph And($k$), that is, $f=1$.  \newline
  (ii) Now, suppose that   $g(W_0)=V_2$. Let $v \in W_0$. We saw that the vertex $-v$ is the unique vertex in  $V_2$ such that its degree in the graph $H$ is equal to the degree of $v$, that is,  $deg_H(-v)$=$deg_H(v)$. It follows that for every vertex $x$ of $H$ we have
$g(x)=-x$. Since $\Gamma$=And($k$) is an Abelian Cayley graph, then the mapping $a:V(\Gamma) \rightarrow V(\Gamma)$ defined by the rule $a(v)=-v$ is an automorphism of the graph $\Gamma$. Let $b=af$. Thus $b$ is an automorphism of $\Gamma$ such that its restriction to $W_0$ is the identity automorphism. Now, by what is proved in (i), it follows that $b=1$. Since $a$ has order 2 in $Aut(\Gamma)$, then $f=a$.\

 We now conclude that if $A=Aut(\Gamma)$ and $A_0$ is the stabilizer subgroup of the vertex $v=0$, then $A_0=\{ 1,a \}$, and hence we have $|A_0|=2$. Now, from Orbit-Stabilizer theorem it follows that $|A|=|A_v||V(\Gamma)|$=$2(3k-1)$.  \

On the other hand, we know that $Aut(\Gamma)$ has a subgroup isomorphic to the cyclic group $\mathbb{Z}_{3k-1}$, that is, $S=\{f_v \ | \ v\in  \mathbb{Z}_{3k-1} \}$, where $f_v:V(\Gamma) \rightarrow V(\Gamma)$, $f_v(x)=x+v$ for every  $x \in \mathbb{Z}_{3k-1}$. It is easy to check that $A_0 \cap S=\{ 1 \}$, hence $| SA_0 |=\frac{|S||A_0|}{|S \cap A_0|}=2(3k-1)=|A|$. Thus,  we have $A=SA_0$. Since the index of $S$ in $A$ is 2, so $S$ is a normal subgroup of $A$. Now we have $A \cong S \rtimes A_0 \cong \mathbb{Z}_{3k-1} \rtimes \mathbb{Z}_2 \cong \mathbb{D}_{2n}$, where $n=3k-1$.

\end{proof}

\end{document}